\documentclass{primus}

\usepackage{amsmath,amssymb}
\usepackage{color}
\usepackage{graphicx}
\usepackage{url}
\usepackage{hyperref}

\title{COLLABORATIVE LEARNING THROUGH FORMATIVE\\ 
PEER REVIEW WITH TECHNOLOGY}

\author{Carrie Diaz Eaton\\
Center for Biodiversity\\
Unity College\\
Unity, Maine 04988, USA\\
ceaton@unity.edu
\and
Stephanie Wade\\
Center for Environmental Arts and Humanities\\
Unity College\\
Unity, Maine 04988, USA\\
swade@unity.edu}


\newcommand{\AmSLaTeX}{$\cal A$\kern-.1667em\lower.5ex\hbox{$\cal
M$}\kern-.125em $\cal S$-\LaTeX}

\newcommand{\matlabs}{{\sc Matlab}{\/ }}
\newcommand{\matlab}{{\sc Matlab}}

\begin{document}


\makePtitlepage
\makePtitle

\begin{abstract}
This paper describes a collaboration between a mathematician and a compositionist who developed a sequence of collaborative writing assignments for calculus. This sequence of developmentally-appropriate assignments presents peer review as a collaborative process that promotes reflection, deepens understanding, and improves exposition. First, we distinguish writing-to-learn from writing-in-the-disciplines. Then, we review collaborative writing pedagogies and explain best practices for teaching peer review.  Finally, we present an implementation plan and examples of student work that illustrate improved understanding of content and improved exposition.
\end{abstract}

\listkeywords
%

\section{INTRODUCTION}
This paper describes a collaboration about collaboration. A mathematics professor and a composition professor devised a series of collaborative writing and peer review activities for a Calculus I and II course sequence. The mathematics professor learned about composition research, which improved her ability to teach writing. The composition professor learned about the disciplinary conventions of writing in applied mathematics, which improved her understanding of writing across the curriculum. She also learned that distinguishing between writing to learn and writing in the disciplines helps faculty stage assignments. The students gained deeper understanding of mathematical content and improved ability to communicate mathematical concepts. We demonstrate how this project connects research in writing pedagogy and mathematics pedagogy.  Then, we present the process and the results of our collaboration. We hope our work will motivate readers to pursue collaborative, transdisciplinary projects. Through collaboration and peer review, we improved our understanding of our work, we improved this paper, and we improved our students' learning.

\section{WRITING TO LEARN VERSUS WRITING IN THE DISCIPLINES}

Collaboration requires clear communication and precise definition of terms. For example, composition studies, a field that focuses on writing and rhetoric,  distinguishes between writing-to-learn (WTL) and writing-in-the-disciplines (WID). Both WTL and WID are aspects of writing across the curriculum (WAC). WTL uses informal writing to promote student engagement and to facilitate learning.  WID teaches the conventions of disciplinary-specific genres \cite{Beidleman95}. Writing-to-learn views students' errors and mistakes as part of the learning process and emphasizes writing as a way to understand content. Examples of WTL activities include journals, informail writing, and online discussion boards or blogs \cite{Cooper12}. Writing-in-the-disciplines builds on writing-to-learn activities, with additional attention to the conventions of disciplinary-specific genres. Typically, WTL activities are most appropriate for introductory classes, while WID assignments are more appropriate for intermediate and advanced classes, after the students have at least basic understanding of their discipline and of the genres expected in their discipline. While WTL and WID are related in that the former can facilitate the latter, conflating the two is problematic, because each has different goals and requires different attention.

Reva Kasman's 2006 article on peer review in mathematics demonstrates the importance of distinguishing between WTL and WID \cite{Kasman06}. Kasman documents an assignment that presented fictional proofs for students to critique. Kasman intentionally created proofs with errors. She reports that finding and correcting these errors improved students' understanding of content, but that students' writing did not improve (6, 11). The activities described in Kasman's study fall into the category of WTL. Helping students improve their writing in addition to their understanding of content requires several further considerations. Research indicates that as students continue to develop their understanding of increasingly complex material, they will make mistakes, which may explain why the writing of the students in Kasman’s study did not improve  \cite{Kroll78, Kutz78}. Furthermore, WID activities should explicitly introduce students to the genre expectations of disciplinary-specific writing, while the assignment in Kasman's study focused on content.

Including both WTL and WID in mathematics classes may initially appear too time-consuming. Mathematics teachers need to cover content. Breaking writing assignments into stages, teaching disciplinary-specific genre expectations, and reading and responding to multiple student drafts all make further demands on in-class and out-of-class time. Collaborative writing and peer review address the problem of time because they present writing as a process while reducing the number of assignments faculty read, evaluate, and grade. Collaborative writing and peer review allow students to measure their own progress, deepen their understanding of content, and improve exposition by providing early feedback and examples of other student work.

\section{COLLABORATION AND PEER REVIEW FOR UNDERGRADUATES}

Collaboration has had a place in writing classes for many years \cite{Bruffee84}.  In the 1970s, educators in the United States began to employ collaborative pedagogies as a practical means of improving student learning, especially for students who resisted traditional pedagogies. These students, some underprepared and some well prepared, responded to peer assistance more so than faculty interventions. 

The history of collaboration is rooted in social constructivist theories of learning \cite{Bruffee84}. Social constructivist theories posit education as a process of making meaning with intellectual and affective dimensions. Collaborative pedagogies align with social constructivist theories because they engage students in active participation, in group work, and in dialogue.\footnote{For a thorough theoretical rationale of the value of collaboration, see Kenneth Bruffee's 1984 article ``Collaborative Learning and the Conversation of Mankind'' \cite{Bruffee84}.} Group writing of assignments engages students in active participation and dialogue, so they gain better understanding of concepts and genre conventions, which prepares them to meet WID expectations.

Group work and group writing assignments are both examples of collaborative pedagogical techniques.  Peer review is another collaborative pedagogical tool that is both efficient and effective. Empirical studies demonstrate the benefits of peer review over feedback from teaching assistants in science classes \cite{Patchan11} and the importance of teaching students to focus on content \cite{Patchan09}.   Reflection on the peer review and composition process facilitates metacognition, which the Council of Writing Program Administrators defines as ``the ability to reflect on one’s own thinking as well as on the individual and cultural processes used to structure knowledge'' \cite{Council of Writing Program Administrators}. The skills of metacognition are important to writing because several studies correlate it with transfer. In other words, when students are more aware of the way they think and act in one writing situation, they become more adept at thinking and writing in other situations \cite{Beaufort07, Downs07}.  

Collaboration between faculty from different disciplines can reinforce transfer. For example, in our collaboration, the composition professor learned more about writing in mathematics, so she could explicitly tell her students how to apply lessons from her class to other situations. The mathematics professor learned about how students learn to write, so she could more effectively build upon what they learned in their first-year communication classes.

Unfortunately, many students and teachers have had negative experiences with peer review and view it as ineffective. In order for peer review to be effective, students need to understand why they are doing it, students need to know how to give good feedback, and students need incentive to take peer review seriously \cite{Henry11}. While many colleges have writing specialists who could provide resources and onsite support, below we offer a brief overview of how to introduce peer review and then we explain how collaborative writing technologies can facilitate successful peer review.

\section{IMPLEMENTING PEER REVIEW FOR COLLABORATIVE LEARNING}

When introducing peer review for the first time, teachers will want to explain to students why they are doing it \cite{Henry11}. For example, a teacher might assign this short video about peer review created by MIT \cite{MITVideo}. This video includes faculty and student voices to convey the three main benefits of peer review: 1. Getting feedback aids revision. 2. Seeing how classmates approached assignments aids revision. 3. Giving feedback aids revision. These three reasons all emphasize how peer review aids revision. Thus, for students to effectively engage in peer review, they need to understand that writing, even writing in mathematics, is a process that involves drafting, revising, and editing. 

Separating revision from editing is one step towards teaching students how to give good feedback. Revision involves choices about content; editing involves choices about language. While editing is important, it is best to save it for the very end of the writing process for two reasons. First, editing the language of a section that might get cut out of the final draft wastes time, and second, focusing on correctness when drafting can create writer's block and prevent writers from attempting to represent complex ideas. If we are employing writing as WTL, we need our students to focus on content first \cite{Sommers82}. 

We also want to address students' possible inexperience with peer review. As Richard Chisholm writes: ``Experienced writers and experienced reviewers know that a solo draft is \textit{only} a draft and that the purpose of peer review is to stimulate the writer to rethink the entire document'' \cite{Chisholm91}. One way to teach inexperienced writers about the relationship between peer review and revision is by asking students to begin peer review by setting an agenda. What this means is that students tell each other where they want feedback and what, if anything, is off limits. Faculty who use a grading rubric might direct students to refer to the rubric when they set their agenda. Setting the agenda not only provides students with some directions for peer review, it also gets students reflecting on their own work, which contributes to metacognition \footnote{Joel A. English \cite{English88}  and L. Lennie Irvin \cite{Irvin04} explain the connections between metacognition, reflection,  and computer-mediated writing.}.

Another way to emphasize content is to tell students to ask questions about their peers' papers. Questions help those who are giving feedback, by offering a way to give feedback without direct criticism. They also help those who are getting feedback by showing them that they have choices in terms of what they say and how they say it. Thus, these questions should address both form and content. For example: ``Why did you put this point in this paragraph?” or``What do you mean by regression?'' Whether using online collaboration tools or working with print copies, students can include these questions as comments on each other's papers. 

In addition to questions, students should be directed to write comments at the end of each other's papers. These comments should take the form of summary, praise, and constructive critique. The summary allows writers to see if their main points are clear and ensures that the reviewers are reading carefully. The praise helps writers see their strengths; it is surprising how often students note that peer review helped them develop confidence in their own writing. Finally, the constructive criticism allows writers to see where they need improvement. Providing students with models of good peer review is another way to make peer review more effective. Chisholm provides such models on pages 16-17 of his paper \cite{Chisholm91}.

While students may complete peer review in class or outside of class, it is useful to devote some time in class for students to set their agenda and to discuss their work.  Allowing the activity to be completed for homework allows students to work at various paces and accommodates slow readers who need more time to compose useful feedback.

The final step in peer review is for students to compose a plan for revision. This takes the form of a paragraph or two or several bullet points about 1) what peer review helped them see about their work; 2) the advice they will follow, the advice they will ignore, and why; and 3) what  next steps they will take to revise. The plan for revision directs students to revise, facilitates their ownership of their writing, and also contributes to metacognition.

In addition, real-time collaborative writing technologies can facilitate collaborative writing and peer review in and out of class. One popular example is GoogleDocs: the suite of editors available within GoogleDrive, \url{https://drive.google.com}, a cloud-based drive for housing files.  GoogleDocs includes such editors as GoogleSpreadsheet (which can run script for programming), GoogleDocument (word processing editor with a \LaTeX-based equation editor), and GooglePresentation. Students can work on every part of a project collaboratively, from computing figures on a spreadsheet, to writing the lab report, to presenting results to the class. The chat window available during the collaboration session allows for editing suggestions, \textit{e.g.} ``Can you do some rephrasing to facilitate that section transition?,'' without compromising the integrity of the working document. Google also offers sharing and permissions options, so users can set very specific controls over who has access to their document and how much access they have. Individuals participating in group reports with full editing privileges have access to detailed revision histories color coded by account user.  This revision history can be useful for both students and teachers when reflecting on group participation or on writing evolution. Peer reviewers can be added with editor or comment-only privileges.   

There are also free, web-based collaborative \LaTeX{ } editors.  Write\LaTeX, \url{https://www.writelatex.com/}, is an example of such an editor that allows for simultaneous cloud editing.  Its split-screen design allows for the \LaTeX{ } editor to appear alongside the PDF document view, and files can be synced with Dropbox for cloud-based access from anywhere. We chose to write the majority of this paper in GoogleDocs, using collaborative writing technology to write about collaborative writing technology.  After our first draft, we inserted our paper into \LaTeX{ } for formatting reasons, then continued the collaborative writing process in Write\LaTeX.  The advantages of writing collaboratively in \LaTeX{ } for mathematical collaborations between math majors and math faculty are clear.  However, GoogleDrive offers a suite of features accessible to for non-math majors who may not need to learn \LaTeX{ } to write collaboratively.

\section{CASE STUDY: COLLABORATION IN A MATHEMATICS SEQUENCE}

At Unity College, all mathematics courses are service courses. Our mathematics courses serve many important objectives: content, technology skills, and, ultimately, the development of a fundamental language to support our sustainability science curriculum. As we strive to meet the needs of our students, we pay close attention to critical thinking and communication skills.
While writing is an important feature of all mathematics courses at \% College, the applied calculus sequence learning outcomes include independently written full laboratory/research reports by the end of the second semester.  Students work on labs in groups, which helps foster learning \cite{Kast93}. Then they write the resulting research reports in Calculus I as a group, and in Calculus II independently.  To ease the transition between semesters, the same grading rubric is used \cite{CBrubric12}; however, in the second semester, the bar for an A is raised, which elevates personal responsibility and requires higher target mastery.  Staging mathematical writing over two semesters follows the recommendations of Kelly and LeDocq's study of writing in mathematics \cite{Kelly01}. 

Through collaborative writing, students in Calculus I work to create a final product that represents the effort and knowledge base of the group.  Through these group writing assignments, they practice the norms of the discipline, which deepens their knowledge of content and of genre conventions. In this way, they practice writing to learn. When they are asked to transition to independent writing in Calculus II, students lose the informal peer group that collaborative writing provides and they are simultaneously asked to face more challenging disciplinary-specific conventions. Some students have no trouble with this transition.  Other students struggle. Formal peer review allows instructors to provide students with formative feedback in the early stages of writing as they make this transition. In this way, students continue to practice writing to learn as a means of improving their ability to write according to disciplinary-specific conventions.

Appendix A includes a worksheet that describes the peer review process. Through class discussion and review of the worksheet, we explain the three benefits of peer review noted above: 1) feedback, 2) improved ability to critique one’s own work, and 3) more familiarity with content and genre.  Students were asked to arrive to class with a draft of their report written. The teacher then randomly created pairs of peer reviewers. (see the next section for variations). Next, the pairs set the agenda for each peer review and shared their papers. Students worked on the peer review for the remainder of class time and  finished it for homework for the next class. The peer review grade included both effort in creating a reasonable first draft for review and effort and thoughtfulness as a reviewer, which can be an additional source of motivation. Revisions of the final draft were due the following week, and the revised product was assessed using the standardized research report rubric developed by the Unity College Center for Biodiversity \cite{CBrubric12}.  

Students in Calculus I who engage in a collaborative project and group writing were asked to reflect on their experience.  The following are some quotations from these students:
\begin{itemize}
	\item``We have been working on another project for calculus, which I really enjoy doing. It's nice to work with peers on something that is applicable in our fields. I think that it's beneficial to have somebody else explain the same thing but in a different way than someone. [Name omitted] has a really great way of simplifying everything and breaking things down so that I actually understand everything that we are doing in our project. I really enjoy having her point of view.''
\item ``The project that we have begun working on this week is very helpful too. Not only am I getting to figure it out myself and use a real-life type scenario, but I am also getting the opportunity to explain the process and everything to someone else. This has helped me to really put what I'm doing into words and make sure both they and I have it down.''
\item ``We used Google docs to write it all up, and its so cool how you can literally work on the project at the same exact time but be in two different places!  Her and I are planning on finishing up the rest of the paper this weekend and getting it in so we don't have to worry about it.''
\end{itemize}
The following are a examples of improved exposition as a result of the peer review process introduced in Calculus II.  To give some context, in our particular lab, we developed a program in GoogleSpreadsheet using script to calculate the area under the graph of a positive function.  A newer, \matlab-based, version of this project is available in Appendix B.
\begin{itemize}
	\item Draft Title: ``Student Programming Project: Trapezoid Rule Project''
	\item Peer Comment: ``This is good as a description of the paper, but parts of it are unhelpful as the `title'. The title should describe what is going on in your paper.''
	\item Final Title: ``Calculating area using Trapezoid Rule Project in Google Documents''
\end{itemize}
The following example is from the introduction of the same paper.  Many students struggle with explaining why they are doing this exploration.
\begin{itemize}
	\item Draft: ``Anti-derivatives are used in integration to find the equation for the integral in question, but there are times in which anti-derivatives are unable to be used.''
	\item Peer Comment: ``This is confusing. Try re-wording.''
	\item Final: ``Anti-derivatives are used in integration to find the equation \textbf{for calculating the area under the curve of the given function}, but there are times in which anti-derivatives are unable to be used.''
\end{itemize}
The writer is still trying to clarify that anti-derivatives cannot always be calculated by hand, but still the exposition of this idea is still much improved.  Even the latter comment, while not particularly as thoughtful, still forces the reviewed student to reflect on what he or she was trying to convey. In both cases, peer review helped the student solidify understanding of the link between using numerical approximations derived using the equivalence of area to antiderivatives.

Sometimes the peer review comments were a much needed early-intervention in the writing process, and therefore led to a much improved student product.  In one case, a student included the code and the function calls and output as an Appendix and only referred to it in the results section.  However, the student had included no explanation for the code, no words that summarized the results.
\begin{itemize}
	\item Peer Comment: ``What are your results?...You need to define your variables and offer more explanation.  Think of it as if you are writing this for an audience that has no idea what your project was to begin with.''
\end{itemize}
The peer reviewer also made a similar comment in the summary of review.  In the interest of conciseness, the entire re-written Results section is not included, however, there was a much improved results section breaking down the results in context to the problem at hand.

\section{ASSESSMENT, REFLECTION, AND REVISION}

Over time, we have made modifications to some of the activities described above.  All our students now have institutional Google Apps for Education accounts and are introduced to GoogleDrive in their freshman seminar course.  Without this extra assistance, it was necessary to spend some class time engaged in creating Google accounts and familiarizing them with the features, though some of this can be abated by using the plethora of tutorials available on the web, see for example \cite{MSUtutorial}.  We now also have \matlabs available so we do not use the GoogleSpreadsheet script.  Although GoogleSpreadsheet has programming capability and collaborative ability, \matlabs is a much easier programming tool to learn and also has more tutorials available online, see for example \cite{Bodinetutorials}.  Since this was the students' first introduction into programming, \matlabs also allowed students to explore more readily.  Although these seem like structural changes unrelated to writing, we acknowledge that writing mathematical ideas is in itself a challenge.  Adding difficult, unfamiliar tools or content mastery requirements to this challenge can sometimes interfere with the students' ability to communicate it \cite{Herrington81}.

In addition to staging levels of writing between courses, we have now incorporated staged writing within each course.  The first lab report in each course is limited to the title, results, discussion, and references sections. In the next lab, students finish the report. This allows students to develop fluency with aspects of the genre as they deepen their content knowledge, using writing to think like mathematicians \cite{Poe11}. Peer review is part of both reports. In the second semester, students have two opportunities before advancing to the final independent project.  During the first iteration of the course, only one peer review per person was assigned, and some students were more helpful than others.  Some students were disappointed by less-than-helpful or late feedback. Now each student is assigned two reports to peer review for each assignment.  This increases the chance of helpful criticism and also gives students the opportunity to elicit a diversity of opinions.

Our collaboration has yielded several important lessons: Let students know that peer review is how academics craft effective writing in all disciplines.  Use the peer review handout provided in the appendix and the MIT video, \cite{MITVideo}, to reinforce this.  Metacognition is important, so provide students with opportunities to reflect on their learning \cite{NRC99, Driscoll11}.  Let students tell their peers that they are depending on them for their feedback.  Peer pressure sometimes works better than instructor pressure. Use the collaborative technology first before assigning the task to students, so you will be prepared to vet their concerns. Each technology has its own peculiarities. In addition, students may find a new technological or pedagogical tool, like GoogleDrive or peer review, more tedious if no one else employs it, so learn what your colleagues are doing and volunteer to teach others about your favorite tools.  

Most importantly, do not limit yourself to colleagues within your department; perhaps your peers in other departments are adapting these tools in their first-year writing curriculum or in other courses. The writers of this paper have benefited tremendously from collaboration and peer review. The problems we face in undergraduate education are complex and interdisciplinary. By being more deliberate in our own learning processes and by sharing this with our students, we explicitly illustrate the deep connections between collaboration, learning, and communication and we model transdisciplinary problem-solving.

\section*{ACKNOWLEDGMENTS}

We would like to acknowledge the work of the Unity College Center for Biodiversity faculty in creating the rubric used to grade biology research reports \cite{CBrubric12}.  We would also like to thank the Gretchen Schaefer, a friend and Instructional Technologist at Husson University and for her patience and help while we learned and implemented new technology in our classes.  We also thank the organizers of this special issue and the MAA Contributed Sessions on Incorporating Writing and Editing into Mathematics Classes at Mathfest for encouraging this work to be submitted and fostering dialogue.

\section*{APPENDICES}
{\bf APPENDIX~A: PEER REVIEW HANDOUT}

Peer Review
What have your past experience with peer review been like? If you have never done it, what do you think you will like/dislike?\\
\vspace{0.5in}
\\What are the benefits of peer review?
\begin{enumerate}
	\item 
	\item 
	\item 
	\item 
\end{enumerate}
Guidelines for successful peer review:
\begin{enumerate}
	\item Set the agenda with your partners. Tell your partners what sort of feedback you want and what sort of feedback you don't want. Listen to the type of feedback your partners want and don’t want.  Write this down. As you read their papers, address their concerns.
	\item As you read, address content and form. Use questions as a way to point out your concerns.
		\begin{itemize}
			\item For example, if you find a section confusing, you might ask: What are you trying to say here?
			\item For example, if the organization is unclear, you might ask: Why did you put the information or pictures in this order?
			\item You will write these questions in the margins of your partners’ paper(s).
		\end{itemize}
	\item Look at most local issues of language in the last draft. You don’t want to wrestle with the language of a section that may get cut from the final copy. If the language interferes with your ability to understand your partners’ papers, let them know via questions.
	\item You should write at least three paragraphs of feedback.  
		\begin{itemize}
			\item Start with summary. This will help your partners see if they conveyed their main points/achieved their purpose. 
			\item Then, praise what is going well. Be specific.
			\item Finally, write about what needs improvement. Be specific, consider your partner’s agenda, and use questions to encourage critical thinking.
		\end{itemize}
	\item After you give and receive feedback, review your own work and create a plan for revision.  In your plan for revision, you will state the next steps you will take to improve your paper. Note which of your peers’ comment you will heed, which seem off base, and what new possibilities you see in your own work.

\end{enumerate}

\newpage
{\bf APPENDIX~B: RESEARCH PROJECT}

\textbf{Student Programming Project: Trapezoid Rule Project}

We have investigated integrals in two ways: by explicitly using antiderivatives and
by approximation using sum of rectangles. For example, we can calculate
$$\int_1^{10}(3x^{-2})dx=3\frac{x^{-1}}{-1}\left|_1^{10} \right.=-3x^{-1}\left|_1^{10}\right.=(-0.3)-(-3)=2.7$$

This integral is easy to calculate since the antiderivative of its integrand, $-3/x$,
can be found exactly. But what if we want to allow the upper limit to be infinity?
$$\int_1^\infty(3x^{-2})dx.$$

Also some important functions do not have antiderivatives that are straightforward to calculate.
Consider this function, which is used in representing normal distributions:
$$p(x)=e^{-x^2/2}/(2 \pi).$$
We do not have the tools to calculate the explicit antiderivative for this function, so we must approximate any finite integral of $p(x)$ numerically.

We will use MatLab to calculate integrals numerically with the trapezoid rule which we developed in class.  The code can be found in the Bodine, Gross, and Lenhart book.
\begin{enumerate}
	\item Estimate using the trapezoid rule the area under the curve of
$$f(x)=3x^{-2}$$
on the interval $[1, 10]$ using 100 subintervals.  You will have to modify “function” in the f.m file in order to do this.  Then run it and put in the appropriate values for a, b and n.

	\item Now calculate the approximate integrals for $b=100, 1000.$  Note: you may have to modify n in order to get good accuracy, so when you state your results, mention what n you used.  You may want to think about what n gives you the same $\Delta x$ as in \#1.

	\item Calculate the exact value of the integrals in \#2 using antiderivatives.  Compare 1-3.  What do you think is the value of $\int_1^{\infty}(3x^{-2})dx$?

	\item Estimate using the trapezoid rule the area under the curve of
$$p(x)=\dfrac{e^{-x^2/2}}{2 \pi}.$$
over the interval [-1, 1] using an appropriate number of subintervals (you choose $n$ based on what gives you enough accuracy).  Hint: $e^x$ is actually exp(x) in \matlab, $\pi$ is pi, and the square root function is sqrt(x).

	\item Now calculate the same integral over [-2,2] for the same number of subintervals.  Then try 10 times as many intervals.

	\item Repeat, but for [-5,5], then [-10,10].  Conjecture the exact value of this integral from $(-\infty, \infty)$.

	\item Write up a research report, with a title, results, and conclusions/discussion.  Also discuss the accuracy of your results as well as comparing and contrasting results.   Are the results what you expected or different?  Include a copy of your modified function files as an appendix.
	\end{enumerate}
\vspace{0.5in}
There will be three phases of research report:
Phase I.  Due Wednesday 2/27.  You must upload your document to Canvas or submit a GoogleDoc url (make sure you change settings to share by link).  

Phase II.  You will be assigned 2 reports to read of peers and give feedback.  Please review by Monday 3/4.

Phase III. Final drafts will be due, uploaded to Canvas by next Monday 3/11.

Follow the peer review guidelines from class.

Your grade will be a combination of completion, effort in draft, peer review quality of comments you make for other people, and the writing quality, insight, and proper conclusions in the final version (all assessed by the Center for Biodiversity rubric at \cite{CBrubric12}).

\section*{BIOGRAPHICAL SKETCHES}

A. Carrie Diaz Eaton received her BA in Mathematics and MA in Interdisciplinary Mathematics from the University of Maine and PhD in Mathematics with a concentration in Mathematical Ecology and Evolutionary Theory from University of Tennessee.  She then returned to Maine and is currently Associate Professor of Mathematics in the Center for Biodiversity at Unity College.  Her goal as a teacher is to get students to see math in the world around them and to use math as a tool to answer scientific questions.  As her twitter profile @mathprofcarrie says, Carrie has an interest in evolution, ecology, social systems modeling, cloud technology, and becoming a better teacher.

\vspace*{.3 true cm} \noindent B. Stephanie Wade earned a BA in psychology from Wesleyan University, a MA in English at the City College of New York, and a PhD in English, with a concentration on composition studies, at Stony Brook University, in New York. Currently director of writing and assistant professor of writing at Unity College, she teaches writing and humanities classes that aim to help students see the choices available to them as writers, thinkers, and citizens.


\begin{thebibliography}{0}

\bibitem{CBrubric12} Baker, A., D. Potter,  A. Remsburg, A. Arnett, E. Batchelder, E. Creaser, C. D. Eaton, E. Latty, A. Phillippi. 2012. Center for Biodiversity research report rubric. Unity College, Unity ME. {\tt \url{http://goo.gl/SDzrq}}. Accessed 14 March 2013.

\bibitem{Bean86} Bean, J. 1986. Engaging Ideas: The Professor's Guide to Integrating Writing, Critical Thinking, and Active Learning in the Classroom.  San Francisco: Jossey-Bass, 1996.  

\bibitem{Beaufort07} Beaufort, A. 2007. \em{College Writing and Beyond: A New Framework for University Writing Instruction.} Logan: Utah State UP, 2007.

\bibitem{Beidleman95} Beidleman, J., D. Jones, P. Well. 1995. Increasing students' conceptual understanding of mathematics through writing. {\em PRIMUS}. 5(4).

\bibitem{Bodinetutorials} Bodine, E. N. Tutorials for \LaTeX{ } and \matlab. {\tt \url{https://sites.google.com/site/profbodine/tutorials}}. Accessed 29 May 2013.

\bibitem{Bruffee84} Bruffee, K. 1984.  Collaborative learning and the conversation of mankind. {\em College English}. 46(7):635-652.

\bibitem{Chisholm91} Chisholm, R. 1991. Introducing students to peer review of writing. {\em WAC Journal}. 3(1):4-19.

\bibitem{Cooper12} Cooper, A. 2012. Today's technologies enhance writing in mathematics. {\em The Clearinghouse: A Journal of Educational Studies.} 85(2).

\bibitem{Council of Writing Program Administrators} Council of Writing Program Administrators, National Association of Teachers of English, National Writing Project. 2011. Framework for success in postsecondary writing. {\tt \url{http://wpacouncil.org/files/framework-for-success-postsecondary-writing.pdf}} Accessed 13 September 2013

\bibitem{Downs07} Douglas, D., E. Wardle. Teaching about writing, righting misconceptions: (Re)Envisioning `first-year composition' as introduction to writing Studies. {\em College Composition and Communication.} 58.4 (2007): 552-84. 

\bibitem{Driscoll11} Driscoll, D. 2011. Connected, disconnected, or uncertain: Student attitudes about future writing contexts and perceptions of transfer from first year writing to the disciplines. {\tt \url{http://wac.colostate.edu/atd/articles/driscoll2011/index.cfm}}. Accessed 16 March 2013.

\bibitem{English88} English, J. 1988. MOO-based metacognition: Incorporating online and offline reflection into the writing process. {\em Kairos}. 3(1). 

\bibitem{MSUtutorial} Getting started with Google Docs, Michigan State University, East Lansing MI. {\tt \url{http://edutech.msu.edu/online/googledocs/googledocs.html}}. Accessed 24 May 2013.

\bibitem{Gere87} Gere, A. R. 1987 Writing groups: History, theory, and implications. Carbondale, IL: Southern Illinois University Press.

\bibitem{Henry11} Henry, L. 2011. Teaching intellectual teamwork in WAC courses through peer review. {\em Currents in Teaching and Learning}. 3(2). {\tt \url{http://www.worcester.edu/Currents/Archives/Volume\_3\_Number\_2/CURRENTSV3N2.pdf}}. Accessed 10 May 2013.

\bibitem{Herrington81} Herrington, A. 1981. Writing to learn: Writing across the disciplines. {\em College English}. 43(4): 379-387.

\bibitem{Irvin04} Irvin,LI. 2004. Reflection in the electronic writing classroom. {\em Computers and Composition Online}. Spring(2004).

\bibitem {Kasman06} Kasman, R. 2006. Critique that! Analytical Writing Assignments in Advanced Mathematics Courses. {\em PRIMUS} 16(1): 1-15.

\bibitem {Kast93} Kast, D. 1993. Collaborative calculus. {\em PRIMUS}. 3(1): 53-61. 

\bibitem {Kelly01} Kelly, S., E. DeLocq, L. Rebecca. 2001. Incorporating writing in an integrated calculus, linear algebra, and differential equations sequence. {\em PRIMUS}. 11(1): 67-87. 

\bibitem{Kroll78} Kroll, B., J. Schafer. 1978. Error analysis and the teaching of composition. {\em College Composition and Communication}. 29(3): 242–248. 

\bibitem{Kutz78} Kutz, E. 1978. Between Students’ Language and Academic Discourse: Interlanguage as Middle Ground. {\em College English}. 48(4):385–396.

\bibitem{McLeod87} McLeod, S. 1987. Some thoughts about feelings: The affective domain and the writing process. {\em CCC}. 38(4).

\bibitem{NRC99} National Research Council. 1999. {\em How people learn: Brain, mind, experience, and school.} Washington, DC: National Academy Press.

\bibitem{MITVideo} No one writes alone: Peer review, a guide for students. Massachusetts Institute of Technology, Cambridge, MA. {\tt \url{http://mit.tv/zXCM6m}}. Accessed 14 March 2013.

\bibitem{Patchan09} Patchan, M., D. Charney, {\em et al.}  2009. A validation study of students’ end comments. {\em The Journal of Writing Research}. 1(2) 124-152.

\bibitem{Patchan11} Patchan, M., C. Schunn. {\em et al.} 2011. Writing in natural sciences: Understanding the effects of different types of reviewers on the writing process. {\em The Journal of Writing Research}. 2(3): 365-393.

\bibitem{Poe11} Poe, M., N. Lerner, J. Craig. 2011. {\em Learning to Communicate in Science and Engineering: Case Studies from MIT.} Cambridge, MA: The MIT Press.

\bibitem{Sommers82} Sommers, N. 1982. {\em College Composition and Communication.} 33(2): 148-156.    

\end{thebibliography}
\end{document}